\documentclass[leqno,12pt]
{article}
\usepackage{amsmath, amsfonts, amssymb}

\newcommand{\QED}{{\flushright{\hfill \rule{2mm}{2mm}}}}
\newcommand\hook{\mathbin{\hbox{\vrule height .5pt width 3.5pt depth 0pt
\vrule height 6pt width .5pt depth 0pt}}}
\def\proclaim #1. #2\par{\medbreak
  \noindent{\bf#1.\enspace}{\sl#2\par}%
  \ifdim\lastskip<\medskipamount \removelastskip\penalty55
        \medskip\fi}
\def\R{\mathbb{R}}

\def\bT{\mathbb{T}}
\def\bR{\mathbb{R}}
\def\bC{\mathbb{C}}
\def\bN{\mathbb {N}}
\def\bP{\mathbb {P}}
\def\bZ{\mathbb {Z}}
\def\cM{\cal M}
\def\cA{{\cal A}}
\def\cD{{\cal D}}

\def\cL{{\cal L}}
\def\Ham{{\rm Ham}}
\def\Stab{{\rm Stab}}
\def\Susp{{\rm Susp}}
\begin{document}

\title{Partially ordered groups and geometry of contact
transformations\thanks{Supported by the
US-Israel Binational Science Foundation grant 94-00302.}
}
\author{Yakov Eliashberg\\Stanford
University \and Leonid Polterovich\\Tel Aviv University}
\maketitle

\begin{abstract}
We prove, for a class of contact manifolds, that
the universal cover of the group of contact
diffeomorphisms carries a
natural partial order.
It leads to a new
viewpoint on geometry and dynamics of contactomorphisms.
It gives rise to invariants of contactomorphisms
which generalize the classical notion of the rotation
number.
Our approach is based on tools of Symplectic Topology.

\end{abstract}

%\medskip

\centerline{{\it 
Dedicated to D.B.~Fuchs on the occasion of his 60th birthday}}

%\medskip

\section{Introduction and main results}

\medskip
\noindent
\subsection*{1.1 Partially ordered groups}

\medskip
\noindent
Let ${\cal D}$ be a group.  A subset
$C\subset {\cal D}$ is called a {\it normal cone} if
\begin{itemize}
\item[(i)] $f\in C, g\in C\Rightarrow fg\in C$
\item[(ii)] $f\in C, h\in {\cal D}\Rightarrow$ $hfh^{-1}$ $\in C$
\item[(iii)] $1\in C$
\end{itemize}
\noindent Given a normal cone $C\subset{\cal D}$, one defines a relation $f\geq g$
 on ${\cal D}$ by
$$
f\geq g\ {\rm if}\ fg^{-1}\in C\enspace.
$$
Clearly this relation is reflexive $(f\geq f)$ and
 transitive $$(f\geq g,g\geq h\Rightarrow f\geq h).$$
   If it is also anti-symmetric $(f\geq g,g\leq f\Rightarrow f=g)$
    then it is a partial order on $\cal D$.  We call it a {\it
     bi-invariant partial order} induced by $C$.
 Notice that the normality of the cone $C$ implies that if $f_1\geq g_1$ and $f_2\geq g_2$
then $f_1f_2\geq g_1g_2$.

\medskip
\noindent Let us describe now  a way to extract numerical invariants from a
bi-invariant partial order on ${\cal D}$.
  An element $f \in C \setminus \{1\}$ is called
{\it a dominant} if for every $g \in {\cal D}$ there
exists a number $p \in \bN$ such that $f^p \geq g$.
For a dominant $f$  and any $g \in {\cal D}$ set
$\gamma_k(f,g) = \inf \{p \in \bZ \;| f^p \geq g^k\}$,
where $k \in \bN$. Notice that
\begin{description}
\item{(i)} the number $\gamma_k=\gamma_k(f,g)$ is finite, and
\item{(ii)} the limit $\gamma(f,g) = \lim_{k \to +\infty} \frac{\gamma_k}{k}$
exists.
 \end{description}
 \noindent Indeed,  choose $q \in \bN$ such that $f^q \geq g^{-1}$.
If $f^p \geq g^k$ then $g^{-k} \geq f^{-p}$, so
$f^{kq} \geq f^{-p}$ and $p \geq -kq$.
Hence $\gamma_k \geq -kq$ and it is finite, which proves (i).
 Since
$f^{\gamma_n}\geq g^n,\ f^{\gamma_m}\geq g^m$ implies
$f^{\gamma_m+\gamma_n}\geq g^{m+n}$ we conclude that
$\gamma_{m+n}\leq \gamma_m+\gamma_n$, i.e. the sequence $\gamma_k$ is subadditive.
Consider now the sequence
$u_k = \gamma_k +kq$ which, as we just showed, is   non-negative.
Clearly
it is also subadditive. This implies existence of the limit
  $\lim_{k\to\infty} \frac{u_k}{k}$,  and hence of the limit (ii).
  \medskip

  \noindent We will call $\gamma(f,g)$  {\it the
relative growth number} of $f$ with respect to $g$. Notice
that the real number $\gamma(f,g)$ can be of any sign, or equal to $0$.
  \medskip

\noindent If   both  $f$ and $g$ are dominants, then $\gamma(g,f)$
%%%
is also defined,
 and the following inequality holds:
 \begin{equation}
\gamma(g,f) \gamma(f,g) \geq 1.
\tag{1.1.A}\label{eqq3}
\end{equation}
\noindent Indeed,   set $\alpha_k=\gamma_k(f,g)$ and $\beta_k=\gamma_k(g,f)$.
Then we have $f^{\alpha_k}\geq g^k$ and $g^{\beta_k}\geq f^k$.
Hence,
$$g^{\alpha_k\beta_k}\geq f^{k\alpha_k}\geq g^{k^2}.$$
Since $g$ is a dominant this implies that
$\alpha_k\beta_k \geq k^2$. Dividing by $k^2$ both parts
of this inequality and passing to the limit when $k \to +\infty$
we get the required inequality 1.1.A.

\medskip

\subsection*{1.2 The universal cover of the group of contactomorphisms}

 \noindent Let $(M,\xi)$ be a
closed connected contact manifold with a {\it co-oriented}
 contact structure.  Let us denote by $\Gamma(M,\xi)$ the identity
component of the group of contactomorphisms of $(M,\xi)$, and by
$\theta:{\cal D}(M,\xi)\to\Gamma(M,\xi)$ the universal cover of
$\Gamma(M,\xi)$ associated with the basepoint 1.  Our starting
observation is that ${\cal D}(M,\xi)$ {\it carries a natural normal
cone.  }
Let $(SM,\omega)$ be the {\it
symplectization } of the contact manifold $(M,\xi)$. Let us
remind that $SM$ is the total space
of a $\R_+$-subbundle of the cotangent bundle $T^*M$, which is
formed by contact forms compatible with
the given co-orientation of $\xi$.
The symplectic
structure $\omega$ on $SM$ is the restriction of the canonical
symplectic form $d(pdq)$ of the cotangent bundle.
  $SM$ also carries a canonical
conformally symplectic $\mathbb{R}_+$-action.
 Every contactomorphism $\varphi\in \Gamma$ uniquely
lifts  to a $\mathbb{R}_+$-equivariant symplectomorphism
$\widetilde\varphi$ of $SM$, and conversely each
$\mathbb{R}_+$-equivariant symplectomorphism  of $SM$ projects to a
contactomorphism of $(M,\xi)$.  A function $F:SM\to\mathbb{R}$ is
called a {\it contact Hamiltonian} if it is homogeneous of degree
1, that is $F(cx)=cF(x)$ for all $c\in \mathbb{R}_+,\ x\in SM$. The
Hamiltonian flow, generated by a
time-dependent contact Hamiltonian is
$\mathbb{R}_+$-equivariant, and thus defines a contact isotopy of
$(M,\xi)$. Any contact isotopy $\{\varphi_t\}$  is generated in
this sense by a uniquely defined time-dependent contact
Hamiltonian $\Phi_t:SM\to\mathbb{R}$. The isotopy $\{\varphi_t\}$ is
called {\it non-negative} if $\Phi_t\geq 0$ for all $t$.
\medskip

\noindent Let us denote  the contact vector field
$\frac{d\varphi_t}{dt}$ by $X_t$.
A   contact form $\alpha$ with $\xi=\{\alpha=0\} $ can be viewed
as a section  $\tilde\alpha:M\to SM$ of the $\R_+$-bundle $SM\to M$,
so that $\tilde\alpha^*(pdq)=\alpha$. A contact Hamiltonian $\Phi_t:
SM\to M$   can be  pulled   back to $M$  by the section $\tilde\alpha$, and
  the  resulting function
$\widetilde\Phi_t=\Phi_t\circ\tilde\alpha_t:M\to\R$ is  also often called a
contact Hamiltonian of $\{\varphi_t\}$ with respect to $\alpha$.
$\widetilde\Phi_t$ can be equivalently defined
by the formula
$$\widetilde\Phi_t(x)=\alpha(X_t(x))\,.$$
In other words,  $\widetilde\Phi_t$  measures
the transversal to $\xi$
component of the contact vector field $X_t$. In particular, the positivity
of the deformation $\varphi_t$ means that the vector field $X_t$
defines the prescribed co-orientation of $\xi$.
Notice that the contact Hamiltonian with respect to $\alpha$,
which identically equals 1 defines the so-called {\it Reeb flow} of
the contact form $\alpha$. The corresponding contact vector field
$X_t$ can be characterized by the equations $X_t\hook d\alpha=0$ and $\alpha(X_t)=1$.
\medskip

\noindent Let    $C(M,\xi)\subset\cal{D}$ be a set
  of those $f\in \cal {D}$, which can be
represented by  a non-negative path    joining  $ 1 $
  with $\theta(f)$.
  It is easy to see that $C(M,\xi)$ is a normal cone in $\cD(M,\xi)$.
   We call it the {\it non-negative normal cone} in $\cD(M,\xi)$.

\bigskip
In the present paper we study the following two problems.

\medskip
\noindent
{\bf Problem 1.2.A.} Let $(M,\xi)$ be a closed
contact manifold with a co-oriented contact structure. Does the
non-negative normal cone
 $C(M,\xi)$ induces a  non-trivial partial order on ${\cal D}(M,\xi)$?

\medskip
\noindent
{\bf Problem 1.2.B.} Calculate or estimate the relative
growth number $\gamma(f,g)$ of a pair of contactomorphisms
in geometric or dynamical
terms.

\medskip
\noindent The first case when the answer to Problem 1.2.A is positive is
provided by the simplest
contact manifold $S^1 = \bR/\bZ$. Its contact structure is
just the field of $0$-dimensional (!) tangent subspaces,
and the co-orientation comes from the orientation of
the circle.
The group ${\cal D} (S^1)$ consists of all
orientation-preserving diffeomorphisms $f : \bR \to \bR$
which satisfy $f(x+1) = f(x)+1$, and the normal cone
$C(S^1)$ is formed by those $f$ which satisfy $f(x) \geq x$ for
all $x \in \bR$. Clearly $C(S^1)$ induces a partial order
on ${\cal D} (S^1)$. Namely, $f \geq g$ provided
$f(x) \geq g(x)$ for all $x$.
For higher-dimensional manifolds Problem 1.2
requires methods of   symplectic topology.
The bridge between this problem and
symplectic topology is given by the following criterion
(see Section 2.1 below for the proof).

\medskip
\noindent
\proclaim Criterion 1.2.C. The relation $\geq$ is a non-trivial
partial order on ${\cal D}(M,\xi)$ if and only if
there are no contractible loops of contactomorphisms
of $(M,\xi)$ generated by a strictly positive
time-periodic contact Hamiltonian.

\medskip
\noindent This criterion can be checked for a class of contact
manifolds. For instance, for the standard contact structure
on $\bR P^{2n+1}$ this is an immediate consequence of
Givental's non-linear Maslov index theory [G] (see 1.3.C below).
In Section 2 we derive the absence of contractible loops
generated by strictly positive time-periodic Hamiltonians
from
the Lagrangian intersection theory
along the lines of [P1, Lemma 3.B]. This enables us to get the positive answer
 to Problem 1.2.A for   spaces of co-oriented contact elements of certain manifolds,
as well as some prequantization spaces (see Section 1.3 below
for precise formulations).
In  Section 1.8  we give  a reformulation of Problem 1.2.A in the language
of symplectic fibrations.
Some    potential
generalizations of our results  are discussed in  Section 1.9.
\medskip

\noindent As far as Problem 1.2.B is concerned, the simplest
case of the circle $S^1$ indicates that the relative
growth can be considered as a contact generalization of the
notion of the rotation number of a diffeomorphism of the
circle (see  Section 1.6 below).
Our main results (see Theorem 1.6.E and  its proof in Section 3.4 below) deal with the
case when $M$ is the space
$\bP_+T^*\bT^n$
  of contact elements of the torus
$\bT^n$. Here we relate
the relative growth
to the stable Gromov-Federer norm and the Mather minimal action.
In order to calculate or estimate the relative growth number
for
$\bP_+T^*\bT^n$ we use the theory of symplectic and contact shapes
developed in [S] and [E1]. This theory provides us with new invariants
of contactomorphisms which turn out to be useful in the study
of our partial order.
The details of this approach are described in Section 3.
Finally, the relative growth gives rise to a canonical
partially ordered {\it metric space} associated with a
contact manifold. This construction is presented in 1.7 below.
 Notice also that in some simple cases (see  Example 1.6.C below)
%%%
the
computation of the invariant $\gamma(f,g)$ is straightforward
 if  the positive answer to Problem 1.2.A is known.

\medskip

\subsection*{1.3 Main examples}

\medskip
\noindent Here we list examples of contact manifolds
$(M,\xi)$ for which we can prove that
 the non-negative normal cone induces the
non-trivial partial order on ${\cal D} (M,\xi)$.

\medskip
\noindent
{\bf 1.3.A. Spaces of co-oriented contact elements.}
Let us
recall that  the {\it space of co-oriented
contact elements}, or, in other
words, the positive projectivization\footnote{Let $E$ be a real
vector space. We say that two vectors $e_1, e_2\in E$ are
equivalent if there exists $  \lambda>0$
 such that $e_1=\lambda e_2$.  The set of all equivalence classes
  (or, in other words, the set of all {\it oriented} lines in $E$)
   is called {\it the positive projectivization of} $E$ and denoted
    by $\mathbb{P}_+E$.  The definition extends in an obvious way
     to vector bundles.}  of a cotangent bundle
     $\mathbb{P}_+T^*Y$ of any smooth manifold $Y$
carries a canonical contact
     structure whose symplectization coincides with the cotangent
     bundle $T^*Y$ with the deleted $0$-section.
The next result is proved in 2.4 below.

\medskip
\noindent
\proclaim Theorem 1.3.B. If a closed manifold $Y$ admits a
non-singular closed 1-form then the non-negative normal cone
induces the non-trivial partial order on
$ {\cal D}(\mathbb{P}_+T^*Y)$.

\medskip
\noindent
{\bf 1.3.C. Prequantization spaces.} Given a closed
symplectic manifold \\ $(W,\Omega)$ with the integral
cohomology class $[\Omega]$, consider a principal
$S^1$-bundle $QW \to W$ whose first Chern class
equals $[\Omega]$. This bundle admits an $S^1$-connection
whose curvature form equals $\Omega$. The distribution
of the horizontal spaces of this connection is an
$S^1$-invariant contact structure on $QW$ transversal
to the fibers. This contact manifold is called a
prequantization space of $(W,\Omega)$. Note that a
given manifold $(W,\Omega)$ may admit different
(in any reasonable category) prequantization spaces.
We refer to [Ki] for the survey on prequantization.

\proclaim Theorem 1.3.D.
Suppose that
$(W,\Omega)$ has a closed Lagrangian submanifold $L$
with the following properties:
\begin{itemize}
\item{} the connection on $QW$ is flat over $L$ (the
Bohr-Sommerfeld condition);
\item{} the relative homotopy group $\pi_2(M,L)$ vanishes.
\end{itemize}
Then the non-negative normal cone
induces the non-trivial partial order on ${\cal D}(M,\xi)$.

\medskip
\noindent For instance such a $L$ exists when $(W,\Omega)$ is the
standard symplectic torus $(\bR^{2n},dp \wedge dq)/\bZ^{2n}$.

\medskip
\noindent
{\bf 1.3.E. Real projective space.}
The standard contact
$\bR P^{2n+1}$ is a particular case of the
prequantization space $QW$. Namely take
$W =  \bC P^{n}$  and let $\Omega$ be the
Fubini-Studi form normalized in such a way that its integral
over the projective line equals 2. Since $\pi_2 (\bC P^n) = \bZ$
this situation is not covered by our previous result.
Nevertheless the relation $\geq$ is a genuine partial order
on ${\cal D} (
\bR P^{2n+1})$. This follows from [G]. In [G] Givental
introduces an   invariant $m(f)$ of an element
$f \in
{\cal D} (
\bR P^{2n+1})$
called the asymptotic non-linear Maslov index. Represent $f$
as the time one map of a Hamiltonian flow $\{f_t\}$ generated
by a time periodic contact Hamiltonian $F$. Intuitively speaking
the number $m(f)$ is defined as the density in $\bR_+$ of the set
of periods of certain closed orbits of the flow
$\{f_t\}$. It is proved in [G] thet $m(f) >0$ provided
$F$ is strictly positive. On the other hand, $m(1) = 0$.
Combining this with 1.2.C above we get the partial
order on
${\cal D} (
\bR P^{2n+1})$.

\medskip
\subsection*{1.4 Dominants}

\noindent
Here we discuss the notion of dominants (see Section 1.1 above)
in the context of contactomorphisms.
In turns out that the group ${\cal D}(M,\xi)$ admits    a natural class of
  dominants.
Let us denote by $C^+ (M,\xi)\subset C(M,\xi)
\subset{\cal D}(M,\xi)$ the set of all elements
which can be generated by a {\it strictly positive} contact
Hamiltonian.

%\begin{proposition}\label{prop2A}
\proclaim Proposition 1.4.A.  Any $f\in C^+ (M,\xi)$
is a dominant.
%\end{proposition}

\medskip
\noindent
This is an immediate consequence of the following
elementary statement.

%\begin{lemma}\label{prop4B}
\proclaim Proposition 1.4.B. Assume that $f,g\in{\cal D}(M,\xi)$.
Then $f\geq g$ if and only if these elements can be generated
 by contact Hamiltonians $F$ and $G$ with $F\geq G$. Moreover
if $f \geq g$ then the Hamiltonians $F$ and $G$ can be chosen
to be time-periodic.
%\end{lemma}

\noindent{\bf Proof:}

\noindent 1) Assume $f\geq g$.
 Then $gf^{-1}$ can be represented by a
 path $h_t$ generated by a non-positive contact
 Hamiltonian $H$.  Suppose that $g_t$ generated by $G$ represents $g$.
  Then set $f_t=h^{-1}_tg_t$.  This path represents $f$ and is generated
   by $F(x,t)=-H(h_tx,t)+G(h_tx,t)$.  Clearly $F\geq G$.  The periodicity
    can be achieved by a suitable time reparametrization near $t=0$ and
     $t=1$.

 \noindent 2)  Assume that $f_t,g_t$ represent $f$ and $g$ and are generated
 by $F$ and $G$ respectively.  Then $g^{-1}_tf_t$ is generated by
  a non-negative Hamiltonian provided
   $F\geq G$.
     \QED
   %\flushright
  % {\hfill \rule{2mm}{2mm}}

\medskip

\subsection*{1.5 Calculation of the relative growth}

\medskip
\noindent
In view of Proposition 1.4.A  for every
$f\in C^+(M,\xi)$ and $g\in {\cal D}(M,\xi)$ one
can define the relative growth
$\gamma(f,g)$. The calculation
of the relative growth seems to be
a non-trivial
problem.
In Section 3 below we present an approach to this problem
in the case when $(M,\xi)$ is the space of co-oriented
contact elements to the $n$-dimensional torus. The approach
is
 based on the
theory of symplectic/contact shape developed in \cite{S},
\cite{E}.
  Here is a sample result. Besides an elementary observation in Example 1.6.C
  this is the only class
of multi-dimensional examples where we can precisely
calculate
the relative growth number.

 \medskip\noindent
Let $(p,q)$ be the canonical
coordinates on $T^*\mathbb{T}^n$, and let $F(p),G(p)$ be two contact
Hamiltonians on $T^*\mathbb{T}^n\setminus\{{\rm zero\
section}\}=S(\mathbb{P}_+T^*\mathbb{T}^n)$.  Assume that $F(p)>0$ for
all $p\neq 0$, and $G(p)$ is strictly positive for some $p\neq 0$.
Then one has (see Section 3.3 below for the proof)

%\begin{theorem}\label{theorem2B}
\proclaim Theorem 1.5.A. Let $f,g\in {\cal D}(M,\xi)$
be elements generated by $F$ and $G$ respectively.
  Then $\gamma(f,g)=\max\limits_{p\neq 0}\ \frac{G(p)}{F(p)}$.
%\end{theorem}

%\medskip

\subsection*{1.6 Relative growth as the generalized rotation number}

\medskip
\noindent
Here we present two specifications of the relative growth
which can be considered as an extension of the classical
notion of the rotation number of a circle diffeomorphism.
Before going into details let us explain the main motivation.
It is a classical dynamical idea to measure the speed of rotation
of the trajectories of a flow around a given cycle in the
manifold.
Such a measurement
can be performed rigorously and
proved to be useful in various situations including
the circle diffeomorphisms and
Hamiltonian dynamics.
In particular it is closely related
to asymptotical properties of the set of periods
of certain closed orbits of the flow.
We take a different point of view and
consider a flow as a curve on the group of diffeomorphisms.
\footnote
{See [P4] for some applications of this viewpoint in the context
of Hofer's geometry.}
Our suggestion is to measure the rotation of this curve
around a cycle in the group!
Using the notion of relative growth we can rigorously implement
this idea   for the group of contactomorphisms.
As we will see   in some examples below, the results
of both measurements (the one we suggest
and the classical one)
are closely related to each other.

\medskip\noindent Here are precise definitions.
Let $\Pi\subset{\cal D}(M,\xi)$ be the (full) lift of
$1\in\Gamma(M,\xi)$, which is identified with the fundamental
group $\pi_1(\Gamma(M,\xi),1)$.  Every element $f\in C^+ (M,\xi)$
gives rise to a function $f \to \gamma (f,e)$ defined on $\Pi$.
Let us mention two properties of this function:
%\begin{equation}
$$f=hg\,h^{-1}\Rightarrow \gamma (f,.) \equiv \gamma (g,.)
\enspace;$$
%\tag{1.5.B}\label{eq3A}
%\end{equation}
%\begin{equation}
$$\gamma (f,e_1e_2)\leq \gamma (f,e_1)+\gamma (f,e_2)
\quad\hbox{for all}\quad
e_1,e_2\in\Pi\,.$$
%\tag{1.5.C}\label{eq3B}
%\end{equation}
\noindent
The first one is obvious.  In order to prove the second
property note that
 the group $\Pi$ is abelian.  Thus
 the inequalities $$f^{\gamma_n(f_,e_1)}\geq e_1^n\ ,\
f^{\gamma_n(f,e_2)}\geq e_2^n$$ imply
$$f^{\gamma_n(f,e_1)+\gamma_n(f,e_2)}\geq (e_1e_2)^n,$$ so that
$$\gamma_n(f,e_1e_2)\leq \gamma_n(f,e_1)+\gamma_n(f,e_2),$$
 and the claim
follows.

\medskip
\noindent Here is  a cousin of the above construction.
Set $\Pi^+ = \Pi \cap C^+(M,\xi)$. Every element
$e \in \Pi^+$ gives rise to a function $f \to \gamma(e,f)$
on ${\cal D} (M,\xi)$. It follows from 1.1.A that
if $f \in C^+$ and $e \in \Pi^+$ then
$\gamma (f,e) \gamma (e,f) \geq 1$.

\medskip
\noindent The next examples 1.6.B--1.6.D clarify the dynamical meaning
of the functions $\gamma (f,e)$ and $\gamma (e,f)$.

\medskip
\noindent
{\bf 1.6.B. Diffeomorphisms of the circle. }
In this case the group $\Pi$ is isomorphic to $\bZ$. Its generator
$e$ is represented
by the loop $x \to x+t, \; t \in [0;1]$. The corresponding
contact Hamiltonian, viewed as a function on
$S^1$ with the contact form $dx$, identically equals 1.
Thus $e \in \Pi^+$.
Denote by ${\rm Rot}(f)$ the rotation number of $f\in
{\cal D} (S^1)$. We claim that
$$\gamma (e,f) = {\rm Rot}(f) = \gamma (f,e)^{-1}.$$
In the first equality $f$ is arbitrary while in the second one
we assume that $f(x) > x$ for all $x \in \bR$.

\medskip\noindent
Let us  prove the second equality. The proof of the first one
is absolutely similar.
First assume that $f^{\gamma_k}\geq e^k$ for some
 $k\in \bN$.  Then $f^{\gamma_k}(x)\geq x+k$ for all $x\in\bR$, so
 ${\rm Rot}(f)\geq \frac k{\gamma_k}$.
  Passing to the limit we get the inequality
$\gamma(f,e)\geq {\rm Rot}(f)^{-1}$.

\noindent Let us verify the opposite inequality.
Suppose that ${\rm Rot}(f)>\frac m\ell$ for some $m,\ell\in\bN$.
We claim that $f^\ell(x)\geq x+m$ for
all $x\in\bR$ (see \cite{CFS}).
Indeed, if this is true for some, but not for all
$x$ then there exists $x_0\in\bR$ such that
$f^\ell(x_0)=x_0+m$.  But then ${\rm Rot}(f)=\frac m\ell$, a
contradiction.  If $f^\ell(x)<x+m$ for all $x\in \bR$
then ${\rm Rot} (f)\leq \frac m\ell$, and again we get
a contradiction.  The claim follows.
 Thus $\gamma_m(f,e)\leq\ell$,
so $\frac{\gamma_m(f,e)}m\leq\frac\ell m$.
Taking now a sequence $\frac m\ell\nearrow
{\rm Rot}(f)$ we get that $\gamma(f,e)\leq\frac 1{{\rm Rot}(f)}$.
This completes the proof.

\medskip
\noindent
{\bf 1.6.C. Reeb flows on prequantization spaces.}
Suppose  that  a contact manifold
$(M,\xi)$ admits a contact form $\alpha$ whose Reeb flow is
$1$-periodic.
  Let $\varphi_t\in{\cal D}$ be the lift to the universal cover of
this flow, and set $e=\varphi_1$.  All prequantization spaces
 defined in 1.3.C  admit a contact form with this property.

If the non-negative normal cone induces a non-trivial
 partial order on ${\cal D}(M,\xi)$ (comp. 1.3.D and 1.3.E above),
then
 {\it  for any $t\in\R$ we have}
$$\gamma(e,\varphi_t)=t\;.$$
Indeed, we have  $\varphi_a\leq\varphi_b$, provided that $a \leq b$,  because
$\varphi_a$ is the time $1$ map of the constant contact Hamiltonian $a$.
On the other hand, if $t=\frac {p}{q}$
then $(\varphi_t)^q=e^p$, and the claim
 follows immediately  for rational numbers,
and for irrational by continuity.

When $M$ is the standard contact $\bR P^{2n+1}$ one can get an
estimate  of $\gamma(f,e)$ and $\gamma(e,
f)$ in terms of the
non-linear Maslov index $m(f)$ (see 1.3.E above). Let
$f \in {\cal D} (\bR P^{2n+1})$ be a lift of
 the time $1$ map of any strictly positive time-independent
Hamiltonian, or, in other words, the time one map of the Reeb flow
of any contact form.
It follows from [G] that
$\gamma (f,e) \geq (n+1)/m(f)$ and
$\gamma (e,f) \geq m(f)/(n+1)$.

\medskip
\noindent
{\bf 1.6.D. Stable norm on $H_1(\bT^n,\bZ)$.}
Note that the torus
$\mathbb{T}^n$ acts on $\mathbb{P}_+T^*\mathbb{T}^n$ by
shifts, thus there
exists a natural monomorphism $\pi_1(\mathbb{T}^n)\hookrightarrow
\pi_1(\Gamma)$. We will identify its image with
$H_1(\mathbb{T}^n,\mathbb{R})$\,).

\medskip\noindent Let $\rho$ be a Riemannian metric on
$\mathbb{T}^n$.
Denote by $\parallel\ \parallel_\rho$ the Gromov-Federer
stable norm on
$H_1(\mathbb{T}^n,\mathbb{R})$ associated to $\rho$.
This norm can be defined in several equivalent  ways
(see sections 4.C,D
and especially $4.20\frac{1}{2}$ in [Gr2]). The simplest
one is as follows.
\footnote
{In 3.4 below we use another definition.}
The stable norm of an integral class
$e \in H_1(\bT^n,\bZ)$ equals
to $\lim_{k \to +{\infty}} \frac {l_k}{k}$,
where $l_k$ is the minimal length of a
closed
geodesic in the class $ke$.

\noindent Let $f \in {\cal D} (\bP_+ T^* \bT^n) $ be the
time-one map of the geodesic flow of $\rho$.

\medskip
\noindent
\proclaim Theorem 1.6.E.
The following inequality holds:
$$
\gamma (f,e)\geq
||e||_\rho$$ for any $ e\in H_1(\bT^n,\bZ)$.

\medskip
\noindent
We refer to Section 3.4 below for the proof and more detailed
discussion.

\medskip
\noindent
\subsection*{1.7 The geometry of the relative growth}

\medskip
\noindent
Consider the following question which,
in fact, has  motivated this paper.
The group of compactly supported Hamiltonian diffeomorphisms
of a symplectic manifold (which can be closed as well as open)
carries a remarkable
geometric structure, the Hofer
biinvariant metric (see [H2],[HZ],[EP],
[MS],[P4]).
There are no
natural biinvariant metrics
on the group of contactomorphisms.
\footnote
{Here is an explanation. Consider the natural embedding
$PSL(2, \bR) \to {\rm Diff} (S^1)$. Every
biinvariant metric on ${\rm Diff (S^1)}$
induces a biinvariant metric on $PSL(2,\bR)$.
But such a metric cannot generate a natural topology on
$PSL(2,\bR)$ since the conjugacy classes are non-compact!}
So it seems reasonable  to ask whether this group
admits a geometric structure at all.
The next construction should be considered as an attempt
to answer this question.

\medskip
\noindent
{\it A partially ordered metric space} is a metric space
$(Z,d)$ endowed with a partial order$\succeq$ such that
for every $a,b,c \in Z$ with $a \succeq b \succeq c$
holds $d(a,c) \geq d(b,c)$.

\medskip
\noindent
Let $(M,\xi)$ be a closed contact manifold.
Assume that $\geq$ is a genuine partial order on
${\cal D} (M,\xi)$. It turns out that in this situation one
can associate with $(M,\xi)$ in a functorial way a partially ordered
metric space $(Z(M,\xi), d,\succeq)$, where the metric $d$
comes from the relative growth $\gamma (f,g)$.
\medskip

\noindent The construction
goes as follows.
\footnote
{In fact one can perform it in a much more general context
of partially ordered groups.}
We will abbreviate $C^+ = C^+ (M,\xi)$  (see 1.4. above) and
${\cal D} = {\cal D} (M,\xi)$.
First, we formulate an elementary lemma, similar to the   inequality 1.1.A above.

\medskip
\noindent
\proclaim Lemma 1.7.A. For every $f,g,h \in C^+$
$$\gamma (f,h) \leq \gamma (f,g) \gamma (g,h).$$

\medskip
\noindent
Let us define a function
$\kappa: C^+ \times C^+ \to [0; +\infty)$
by
$$\kappa (f,g) = \max (\log \gamma (f,g), \log \gamma (g,f)).$$
Obviously $\kappa$ is symmetric and vanishes on the diagonal.
Further,
it follows from the inequality 1.1.A that $\kappa$ is non-negative,
and from
Lemma 1.7.A that $\kappa$ satisfies the triangle inequality.
Thus $\kappa$ is a pseudo-distance. We say that two
elements $f$ and $g$ in $C^+$ are equivalent if
$\kappa (f,g) = 0$, and consider the  corresponding
quotient space $Z$ of $C^+$.   The pseudo-distance $\kappa$ projects to
a genuine distance function, denoted $d$, on $Z$.

\medskip\noindent Let us define now a
partial order $\succeq$ on $Z$ as follows.
We write $[f] \succeq [g]$ if $\gamma (f,g) \leq 1$.
Then   $(Z,d,\succeq)$ is a partially ordered metric
space, which easily follows from Lemma   1.7.A.
  This completes our construction.

\medskip
\noindent
Let us make two remarks. First, the  natural
projection $C^+ \to Z,\; f \to [f]$  is clearly
monotone with respect to the partial orders on $C^+$ and $Z$.
Moreover,
{\it it is constant on the conjugacy classes
in $C^+$.}  Indeed, assume that $g = hfh^{-1}$
where $f,g \in C^+$ and $h \in {\cal D}$. Take
a positive integer $k$ such that $f^k \geq h$ and $f^k \geq h^{-1}$.
Thus for all natural numbers $n>2k$ one has
$f^{n+2k} \geq g^n \geq f^{n-2k}$. Therefore,
$\gamma (f,g) = \gamma (g,f) = 1$, so $[f] =[g]$.
\medskip

\noindent
Second, the space $Z$
 admits an action of the multiplicative
semigroup of positive integers by isomorphisms, i.e. order
preserving isometries. For $m \in \bN$ define a map
$\Phi_m : Z \to Z$ by $[f] \to [f^m]$. It is straightforward to check that
this map is isometric and monotone with respect to $\succeq$.
 Of course, $\Phi_m \circ \Phi_n = \Phi_{mn}$
for all $m,n \in \bN$.
\medskip

\noindent The geometry of the space $Z(M,\xi)$
is yet to be explored.
\footnote
{For instance,
is there a meaningful description of geodesics on $Z$?}
Let us give two examples.

\medskip
\noindent
{\bf Example 1.7.B. The circle.}
We claim that $(Z(S^1), d,\succeq)$ is simply the Euclidean
line $\bR^1$ endowed with the natural order.
The natural projection $C^+ \to Z$ is given by
$f \to \log {\rm Rot} (f)$. This follows from the fact
that for all $f,g \in C^+ $
\begin{equation}
\gamma (f,g) = \frac{{\rm Rot}(g)}{{\rm Rot}(f)}.
\tag{1.7.C}\label{eqq1}
\end{equation}
The proof of   (1.7.C) goes as follows. Recall from
Example 1.6.B that there exists an element $e \in {\cal D}$
which satisfies
\begin{equation}
\gamma (e,f) \gamma (f,e) = 1
\tag{1.7.D}\label{eqq2}
\end{equation}
for all $f \in C^+$.
Applying Lemma 1.7.A we get that
$$\gamma (f,g) \leq \gamma (f,e) \gamma (e,g) =
\frac{\gamma(e,g)}{\gamma(e,f)},$$
and similarly
$$\gamma (g,f) \leq \frac{\gamma(e,f)}{\gamma (e,g)}.$$
Multiplying these two inequalities and comparing the result
with 1.1.A above we get that in fact each of them must
be an equality! Recall now from  Example 1.6.B that $\gamma (e,f) =
{\rm Rot}(f)$. This proves (1.7.C), and hence confirms our
description of $Z(S^1)$.

\medskip
\noindent Notice that the isometry $\Phi_n$ in this case  acts on $\R=Z(S^1)$ as
the translation $x\mapsto x+\log n$, $x\in\R$.
  In fact, our proof shows that if for a contact manifold
$(M,\xi)$ there exists an element $e \in C^+$ which satisfies
(1.7.D) for all $f \in C^+$ then $Z(M,\xi)$ can be identified
with a subset of the  Euclidean line, invariant under translations
$\Phi_n$, and the natural projection
is given by $[f] = \log \gamma (e,f)$. We will see in the next example
that no such element exists when $M = \bP_+T^*\bT^n$ with $n > 1$.

\medskip
\noindent
{\bf Example 1.7.E. $
\bP_+T^*\bT^n$.}
Consider the subset $K \subset {\cal D}$ which consists
of the natural lifts of the time-one-maps $\phi_F$ of autonomous
contact flows generated by strictly positive contact Hamiltonians
$F = F(p)$. Intuitively speaking this is the positive part
of ``the maximal torus" of the group $\cal D$.
Denote by $Z_K$ the image of $K$ in
$Z(
\bP_+T^*\bT^n )$.

\medskip
\noindent
\proclaim Theorem 1.7.F.
The space $Z_K$, endowed with
the induced metric and partial order, is isomorphic
(in the category of partially ordered metric spaces)
to the linear space $C^{\infty}(S^{n-1})$ of all smooth functions on
the sphere $S^{n-1}$ with the norm $||u||= \max |u|$ and
with the natural partial order.
\footnote
{The natural partial order on $C^{\infty}(S^{n-1})$
is defined as follows.
We write
 $u \succeq v$ for functions $u$ and $v$ when
$u(x) \geq v(x)$ for all $x$. }

\medskip
\noindent
The isomorphism, however, is not canonical.

\medskip
\noindent
{\bf Proof:}
Take an arbitrary strictly positive contact Hamiltonian
$H(p)$. Consider a map
$j: Z_K \to C^{\infty}(S^{n-1})$  given by
$j([\phi_F]) = \log (F/H).$
It view of Theorem 1.5.A
$\gamma (\phi_F,\phi_G) = \max_{p\neq 0} \frac{G(p)}{F(p)}.$
Therefore $j$ is a monotone
isometric bijection.
\QED

\noindent
The conclusion is that the space
$Z (\bP_+T^*\bT^n)$
contains
an infinite-dimensional flat piece.

\medskip
\noindent
\subsection*{1.8 Fat symplectic connections
\footnote
{We refer to [GLS] and [MS] for basic theory of symplectic fibrations.}}

\medskip
\noindent
It turns out that the existence of
the partial order on the universal cover
of the contactomorphisms group is closely related to
the notion of fat connection on a symplectic fibration.

Let $p:P \to S^2$ be a symplectic fibration with the fiber
$(SM,\omega)$
whose structural group is the group of $\bR_+$-equivariant
Hamiltonian diffeomorphisms of $SM$ (and thus the structural group
can be identified with $\Gamma (M,\xi)$).
We call it {\it an $\bR_+$-equivariant symplectic fibration.}
The fibration
$P$ is endowed with the canonical fiberwise symplectic
structure $\omega_x,\; x \in S^2,$ and the
canonical $\bR_+$-action $(x,z) \to (x,R_c z)$ where
$x \in S^2$, $z \in p^{-1}(x)$ and $c \in \bR_+$.
Let $\nu$ be a connection on $P$ whose parallel
transport maps act by equivariant symplectomorphisms.
We call $\nu$ an equivariant symplectic connection
on $P$. The curvature $\rho$ of $\nu$
is a 2-form on $S^2$ which at a point $x \in S^2$
 takes values in the Lie
algebra of the group of $\bR_+$-equivariant
Hamiltonian diffeomorphisms of $p^{-1}(x)$.
In other words, for $v,w \in T_xS^2$ one considers
$\rho (v,w)$ as a contact
Hamiltonian on $p^{-1}(x)$. Consider the splitting
$$TP = T(SM) \oplus TS^2$$
associated to connection $\nu$. Define a 2-form
$\delta$ on $P$ by
$$\delta = \omega_x \oplus -\rho .$$
The form $\delta$ is called {\it the coupling form} of the
connection $\nu$. It is known to be closed. Moreover
$R_c^* \delta = c\delta$ for all $c \in \bR_+$, and
thus $\delta$ is exact. Following Weinstein [W]
we call
the connection $\nu$ {\bf fat} if the coupling
form $\delta$ is symplectic. Note that if $P$
admits a fat connection then the quotient space $P/ \bR _+$
admits a contact structure which extends
the contact structure on fibers $p^{-1}(x)/\bR_+$.

Loops of Hamiltonian
diffeomorphisms are closely related
to
symplectic fibrations over $S^2$
(see [Se],[P3],[P4],[LMP],[M1]). This link admits a
straightforward generalization
to the category of $\bR_+$-equivariant symplectic fibrations.
In particular, there exists a one-to-one correspondence
between $\pi_1 (\Gamma (M,\xi))$ and isomorphism classes
of $\bR_+$-equivariant symplectic fibrations. Denote
by $P(\alpha)$ the symplectic fibration corresponding to
an element $\alpha \in \pi_1 (\Gamma (M,\xi))$. Note that
the fibration
corresponding to the neutral element of
$\pi_1 (\Gamma (M,\xi))$ is simply the trivial fibration
$SM \times S^2 \to S^2$.
Further,
one can show along the lines of [P3] that
an element $\alpha$ can be represented by a loop
generated by a strictly positive contact Hamiltonian
if and only if fibration $P(\alpha)$ admits a fat connection.
Combining this with Criterion 1.2.C we get the following result.

\medskip
\noindent
\proclaim Proposition 1.8.A.
The relation $\geq$ on ${\cal D}(M,\xi)$ is a genuine partial
order if an only
if the trivial
$\bR_+$-equivariant symplectic
fibration $SM \times S^2 \to S^2$ does not admit a fat connection.

\medskip
\noindent
\subsection*{1.9 Possible generalizations}

It is possible that Problem 1.2.A has the positive answer
for all contact manifolds,
but at the moment the authors do not see any tools to handle the problem
in such generality. In fact, overtwisted contact structures on $3$-manifolds
could potentially provide counter-examples to such a general conjecture.

 \medskip\noindent
However, it seems that the contact homology theory, which is currently
under construction by A. Givental, H. Hofer and one of the authors (see [E2]
and [U]) would
provide an adequate tool for  establishing existence of  the partial order
on a large class of those contact manifolds for which
the contact homology algebra is non-trivial.

\medskip\noindent
The simplest closed contact manifold   which
potentially can be covered by this techniques, but not covered
by the results of this paper,
is the standard contact
3-sphere $S^3$.
Let us elaborate on how the contact homology theory could be
applied to Problem 1.2.A.

\medskip\noindent
With   each
element $f  \in {\cal D}(M,\xi=\{\alpha=0\})$ one can associate
(see Section 3 below) a domain $V^+(f)$
in the contact manifold \footnote{
This contact manifold
is
called the stabilization of $(M,\xi)$, see Section 2.2 below.}
$$ (M \times T^*S^1,\tilde\xi=\{\alpha+rdt=0\})$$
 defined
canonically up to a
contact isotopy. An important property of this correspondence
is that if $f \geq g$ then there exists a  contact
  isotopy which takes $V^+(g)$ inside $V^+(f)$. Thus
one needs a technique which
provides non-squeezing type results for domains of the form
$V^+(f)$.
In symplectic topology results of this type can be proved via
the theory of
symplectic capacities.
\footnote
{A symplectic capacity $c$ is an invariant of a symplectic
domain $(V,\Omega)$ which is monotone with respect to inclusion
and satisfies $c(V,\lambda \Omega) = |\lambda| c(V, \Omega)$
(see [HZ],[MS]).} On the other hand,
the only  so far known invariants of contact domains    which  are  suitable for
this job are
 the, so-called, contact shapes (see
 [E1] and Section 3 below). We employ
these invariants in the current paper.
Since the contact
homology provide  more powerful
and adequate invariants of contact domains,
we hope that they
  should  allow us to settle Problem 1.2.A for
 a  more general class of contact manifolds
which includes the standard contact spheres.

 \medskip
\noindent It is also possible  that the finite-dimensional approach
of the theory of generating functions (see [EG]),
employed by Givental [G] in his theory of the non-linear Maslov index
can be extended to a larger class of contact manifolds.
Let us   mention in this context,  that Mohan Bhupal informed us
that he proved, using generating functions,
 existence of the partial order on the group ${\cal D}_0(\bR^{2n+1})$
 of compactly supported contact
transformations of  the standard
contact   $\bR^{2n+1}$.

\noindent Further, one can try to attack Problem 1.2.A using
its reformulation
in the language of symplectic fibrations (see  Section 1.8 above).
In the theory of {\it compact} symplectic fibrations over $S^2$
there exists a powerful technique of Gromov-Witten invariants
(see [Se],[LMP],[M1],[M2])
which leads to various geometric and topological
consequences. It is a challenging problem to extend it to
the $\bR_+$-equivariant case and to apply   to
the study of the partial order on $\cD (M,\xi)$.

\medskip\noindent
Finally, one can consider Problem 1.2.A in a
more general context of Legendrian submanifolds.
The binary relation $\geq$ on the group $\cD (M,\xi)$
admits a natural extension to the following homogeneous
space of the group. Let $\cL$ be a connected component
of the space of all Legendrian submanifolds of $M$.
The tangent space to $\cL$ at some point $L \in \cL$
can be canonically identified with the space
$C^{\infty}(L)$ of smooth functions on $L$.
Consider the field of tangent cones to $\cL$
formed by non-negative functions. This field
is invariant under the action of the contactomorphism
group. Its lift to the universal cover $\tilde{\cL}$
defines in the obvious way a binary relation on
$\tilde{\cL}$. It
would be interesting to decide when
this is a partial order relation.

\medskip
\noindent
\subsection*{1.10. Historical remarks.}

\medskip
\noindent
The notion of a partially ordered group is
classical (see for instance [F]). The relative growth number
$\gamma (f,g)$ is a variation on the theme of {\it order-unit
norm} in partially ordered Abelian groups (see [Go]).
\footnote
{We were unable to find our $\gamma (f,g)$ in the literature.}
For results and references on
partially ordered finite dimensional Lie groups and homogeneous
spaces we refer to [HiN].
The  partial order associated to
 the cone of non-negative Hamiltonians in the group $\Ham (\bR^{2n})$
of compactly supported Hamiltonian
diffeomorphisms of $\bR^{2n}$ was first described without a proof in [E3],
but the first  published proofs appeared in Viterbo's paper [V] and
Hofer-Zehnder's book [HZ].
 The both proofs are based on a construction of a remarkable invariant
\footnote{ In fact it is still an open problem
to show that the
invariants  constructed in [HZ] and [V] coincide.}
of a Hamiltonian diffeomorphism $f$. It equals
the symplectic action of a specially chosen critical point
of $f$. The choice of the critical point is performed
on the basis of a careful study of either the Floer homology
of the action functional of $f$
or (in the similar vein) the Morse theoretical properties
of the generating function of $f$.
A crucial property
of this invariant is that it vanishes provided $f \leq 1$
and is strictly positive provided $f > 1$.
It would be  interesting  to extend the
partial order on $\Ham (\bR^{2n})$
to other
{\it open} symplectic manifolds. Notice that
for {\it closed} symplectic manifolds
no natural partial order on $\Ham$ is known.
The relative growth number $\gamma(f,g)$, and the  geometry
and dynamics related to it,  were
 not yet studied   in the symplectic context,
 %%%
although it seems  that
the technique developed in [HZ] and [V] should allow
 to compute or estimate
  $\gamma (f,g)$
for suitable pairs of Hamiltonian diffeomorphisms.

\section{Establishing the partial order}

 \subsection*{2.1. Loops of contactomorphisms}

  \noindent A {\it loop of contactomorphisms} is a smooth
map $S^1=\bR/\bZ\to\Gamma(M,\xi)$ which takes
  0 to $1\in\Gamma$.  Loops are generated by 1-periodic
time dependent contact Hamiltonians.
  A loop is called {\it non-negative} if
its contact Hamiltonian is non-negative, and
  {\it strictly positive} if its contact
  Hamiltonian is strictly positive.
We start with the following elementary

\medskip
\noindent
\proclaim Proposition 2.1.A. The relation $\geq$ on
${\cal D}(M,\xi)$ is a non-trivial partial order
if and only if every non-negative contractible loop of contactomorphisms
is the constant loop.

\medskip
\noindent
{\bf Proof:}
Assume that every contractible non-negative loop
of contactomorphisms is constant.
Consider an element
$f \in {\cal D}(M,\xi)$ such that $f \geq 1$ and $f \leq 1$.
Then there exist two paths
$\{f^\prime_t\},\{f^{\prime\prime}_t\},
 t\in [0;1]$ which represent $f$ with the following properties.
First,
 $f^\prime_0=f^{\prime\prime}_0=1,\ f^\prime_1=f^{\prime\prime}_1$
  and the paths are homotopic with fixed
end points. Second,
   $\{f^\prime _t\}$ is  generated by a non-negative contact
 Hamiltonian,
and $\{f^{\prime\prime}_t\}$ is generated by a non-positive one.
   We have to show that then both paths are constant.
Indeed, this would mean that
    $f\geq 1$ and $f\leq 1$ implies $f=1$, that
is the relation $\geq$ is antisymmetric.
Without loss of generality one can
assume that $f^\prime_t=f^{\prime\prime}_t=1$ when $t$ is close to 0,
 and $f^\prime_t=f^{\prime\prime}_t=f^\prime _1$
when $t$ is close to 1.  Consider the union of $\{f^\prime_t\}$
 with $\{f^{\prime\prime}_t\}$, where the
second path is taken with the opposite
 orientation.
 We get a loop of contactomorphisms which
is both contractible and non-negative.
By our assumption such a loop must be constant.
Thus $f = 1$ and our claim follows.
The proof of the converse statement is analogous.
\QED

\medskip
\noindent
%\begin{proposition}\label{prop9D}
\proclaim Proposition 2.1.B. If a closed contact manifold
$(M,\xi)$ admits a non-constant contractible
non-negative loop of contactomorphisms, then it admits a contractible
 strictly positive loop of contactomorphisms.
%\end{proposition}

\medskip
\noindent
{\bf Proof:}
The proof is based on ``ergodic" ideas from \cite{P2}.
Let $f:S^1\to\Gamma$ be a smooth contractible
loop, and $F(x,t)$ is the corresponding non-negative Hamiltonian.
 Since the loop is non-constant there
 exists $t_0$ such that $F(x,t_0)\not\equiv 0$.
 In what follows we describe a
 three-step sequence of modifications
of $\{f_t\}$ to a strictly positive contractible loop.

 \medskip\noindent
1)  Set $g(t)=f(t+t_0) f(t_0)^{-1}$.
 The loop $g$ is generated by the Hamiltonian $G(x,t)=F(x,t+t_0)$.
 Hence $G(x,0)\not\equiv 0$.
 We assume that $G(x,0)>0$ for all $x\in SU$,
where $SU$ is the symplectization of
  a domain $U\subset M$.
\medskip

\noindent 2) Let $\varphi_1,\ldots,\varphi_d$ be a
sequence of elements of $\Gamma$ which will be chosen later on.
 Set
$$
h(t)=g(t)\varphi_1g(t)\cdot\ldots\cdot
\varphi_dg(t)\ (\varphi_1\ldots\varphi_d)^{-1}\enspace.
$$
This is a contractible loop generated by the Hamiltonian
$$
H(x,t)=G(x,t)+G(\widetilde\varphi^{-1}_1\widetilde{g(t)}^{-1}x,t)
+\ldots
+G(\widetilde\varphi^{-1}_d\widetilde{g(t)}^{-1}\ldots
\widetilde\varphi^{-1}_1\widetilde{g(t)}^{-1}x,t),
$$
where we write $\widetilde \phi$ for the symplectization
of a contactomorphism $\phi$.
Set $\psi_0=1,\ \psi_k=\varphi_1\ldots \varphi_k$ where $k=1,\ldots,d$.
 Thus $H(x,0)=\sum\limits^{d}_{k=0} G(\widetilde\psi^{-1}_kx,0)$.
 Assume now that $\{\varphi_k\}$ is chosen in such a way that
  $\bigcup\limits^{d}_{k=0} \psi_k(U)=M$.
   Then for every $x\in M$ there exists $k$ such
that $\psi^{-1}_kx\in U$.  Thus $H(x,0)>0$ for all $x\in SM$.
\medskip

\noindent 3) The last inequality
implies that there exists an open arc $\Delta\subset S^1,\ \Delta\ni 0$
 such that $H(x,t)>0$ for
all $x\in SM,t\in\Delta$.
Choose a real number $\alpha$ and a positive integer
 $m$ such that the sequence $\{t+k\alpha\},\ k=0,\ldots,m$
meets $\Delta$ for all $t\in S^1$.  Consider a new loop
$$
e(t)=h(t)h(t+\alpha)\ldots h(t+m\alpha)\
(h(0)h(\alpha)\ldots h(m\alpha)\,)^{-1}\enspace.
$$
Clearly it is contractible.  Its Hamiltonian $E$ is given by
$$
\begin{array}{rr}
& H(x,t)+H(\widetilde{h(t)}^{-1},t+\alpha)+
H(\widetilde{h(t+\alpha)}^{-1}\widetilde{h(t)}^{-1}x,t+2\alpha)+\\
& \ldots +H(h\widetilde{(t+(m-1)\alpha)}^{-1}
\cdot\ldots\cdot \widetilde{h(t)}^{-1}x, t+m\alpha).
\end{array}
$$
Our choice of $\alpha$ guarantees that for all $x\in SM,\ t
\in S^1$ at least one of the
summands is strictly positive.  Since all
other summands are non-negative,
we conclude that $E(x,t)>0$ for all $x,t$.
 This completes the proof.
 \QED

\medskip
\noindent
Criterion 1.2.C follows immediately from 2.1.A and 2.1.B.

\medskip
\subsection*{2.2 Stabilization in the contact category}

\medskip
\noindent
We start with the following well known remark. Let
$(X,\Omega)$ be a symplectic manifold endowed with
a free $\bR_+$-action $x \to R_c x$ ($c \in \bR_+$)
which admits a global slice and
such that $R_c^* \Omega = c\Omega$ for every $c$.
Then the quotient space $X/\bR_+$ carries in a
canonical way a contact structure $\eta
$ such that
$S(X/\bR_+,\eta) = (X,\Omega)$.

\medskip\noindent
Let $(M,\xi)$ be a contact manifold with a co-oriented contact
structure. We   define its stabilization with respect to
dimension as
$$\Stab M = (SM \times T^*S^1, \omega + dr \wedge dt)/\bR_+,$$
where the action of $\bR_+$ is the diagonal one. This action
is given by
$$(x,r,t) \to (R_cx,cr,t)$$
for $c \in \bR_+$.

\medskip\noindent
If the contact structure $\xi$ is defined by a contact $1$-form $\alpha$
then the contact manifold $\Stab M$ is contactomorphic
to $M\times T^*S^1$ with the contact structure defined by the form
$\alpha+rdt$.
\medskip

\medskip\noindent
Let $\tau: SM \to M$ and $\sigma: SM \times T^*S^1 \to \Stab M$
be the natural projections. Denote by $Z$ the zero section
$\{r=0\}$ of $T^*S^1$.
The procedure of stabilization extends naturally to two remarkable
classes of submanifolds of $(M,\xi)$, {\it pre-Lagrangian}
and {\it Legendrian}.
Recall that $K \subset M$ is called pre-Lagrangian if
there exists a Lagrangian submanifold $L \subset SM$ such that
$K = \tau (L)$ and $\tau {\big |}_L :L \to K$ is a diffeomorphism.
The submanifold $L$ is called a {\it Lagrangian lift} of $K$.
Define the stabilization
$$\Stab K = \sigma (L \times Z).$$
This definition does not depend on the choice of
the particular lift $L$ of $K$. Clearly,
$\Stab K$ is a pre-Lagrangian submanifold of $\Stab M$.
\medskip

\noindent
A Legendrian submanifold $\Lambda \subset M$, i.e. a maximal integral
submanifold of $\xi$, can be characterized  by the property
that its lift  $S\Lambda = \tau^{-1} \Lambda$ to $SM$
is Lagrangian.
The Legendrian submanifold
$$\Stab \Lambda= \sigma (S\Lambda \times Z)\subset \Stab M$$
is called the stabilization of $\Lambda$.

\medskip\noindent
Stabilizations arise naturally in the study of loops of
contactomorphisms. Let $N \subset SM$ be a Lagrangian
submanifold, and let $\{h_t\}$ be a loop of contactomorphisms
of $(M,\xi)$ generated by a contact Hamiltonian $H(x,t)$.
Define the suspension $\Susp_h N$ as the image of the
Lagrangian embedding
$$N \times S^1 \to SM \times T^*S^1,$$
$$(x,t) \to ({\tilde h}_tx, -H({\tilde h}_tx,t),t),$$
where ${\tilde h}_t$ stands for the lift of $h_t$ to a
$\bR_+$-equivariant Hamiltonian diffeomorphism of $SM$.
Let $K$ be a pre-Lagrangian submanifold of $M$.
Set $\Susp_h K = \sigma (\Susp _h L)$, where $L$ is a Lagrangian
lift of $K$.  The submanifold $\Susp_h K$ does not depend
on the particular choice of lift $L$, and is a pre-Lagrangian
submanifold of $\Stab M$. Similarly, let $\Lambda$ be a Legendrian
submanifold of $M$. Set
$\Susp_h \Lambda = \sigma( \Susp_h S\Lambda)$. Again, $\Susp_h \Lambda$
is a Legendrian submanifold of $\Stab M$.
Note that if $h$ is the constant loop $h_t \equiv 1$ then
$\Susp_h K = \Stab K$ and $\Susp_h \Lambda = \Stab \Lambda$.

\medskip

\subsection*{2.3 Stable intersections in contact manifolds}

\medskip
\noindent
Let $(M,\xi)$ be a contact manifold (not necessarily compact),
and
$K\subset M$ be a closed
pre-Lagrangian submanifold.
Let $A \subset M$ be a closed submanifold which is either
pre-Lagrangian or Legendrian.
We say that the pair $(K,A)$ has {\it the intersection property}
if
for every contactomorphism $\phi \in \Gamma (M,\xi)$ the intersection
$\phi(K) \cap A$ is non-empty. We say that $(K,A)$
has {\it the stable intersection property} if
$(\Stab K,\Stab A)$ has the intersection property
in $\Stab M$.

\medskip
\noindent The detection and proof of the intersection property is
one of the central problems in
symplectic and
 contact geometry and certain techniques
(Floer homology, generating functions etc.) are developed
to handle it. If it is possible to prove
the intersection property for $(K,A)$ then usually
(but not always) the same argument allows to handle the
stable intersection problem. In 2.4 below we discuss
some examples of pairs with stable intersection
property.

The next result links stable intersections with the partial
order.

\medskip
\noindent
\proclaim Theorem 2.3.A. Suppose that
$(M,\xi)$ contains a pair with the stable intersection
property. Then $\geq$ is a non-trivial partial order
on ${\cal D} (M,\xi)$.

\medskip
\noindent
This is an immediate consequence of Criterion 1.2.C
and the next

\medskip
\noindent
\proclaim Proposition 2.3.B. Let $(K,A)$
be a pair with the stable intersection property.
Let $H(x,t)$ be a 1-periodic Hamiltonian generating a
contractible loop of contactomorphisms. Then there exists a point
$(x_0,t_0) \in A \times S^1$ such that
the function $H(.,t_0)$ vanishes on the ray $\tau^{-1}(x_0)$.

\medskip
\noindent
This statement and its proof is analogous to [P1, Lemma 3.B].

\medskip
\noindent
{\bf Proof of 2.3.B:} Denote by $h$ the loop
of contactomorphisms generated by $H$. Since $h$ is contractible
then $\Susp_h K$ is isotopic to $\Susp_1 K = \Stab K$
through pre-Lagrangian submanifolds of $\Stab M$.
One can easily check that this isotopy extends to a contact isotopy
of $\Stab M$. Thus the stable intersection property guarantees
that $\Susp_h K \cap \Stab A$ is non-empty. Using
definitions of the stabilization and suspension we obtain that
there exists a point $(y_0,t_0) \in K \times S^1$
such that $h_{t_0}y_0 \in A$ and $H(.,t_0)$ vanishes
on the ray $\tau^{-1} (h_{t_0}y_0)$.
Thus the points
$x_0 = h_{t_0}y_0$ and $t_0$ are as needed.
This completes the proof.
\QED

\medskip
\subsection*{2.4. Examples of stable intersections}

\medskip
\noindent
Examples 2.4.A and 2.4.B below combined with Theorem 2.3.B
prove Theorems 1.3.B and 1.3.D, respectively.

\medskip
\noindent
{\bf Example 2.4.A.} Let $Y$ be a closed manifold which admits
a non-singular closed 1-form. The graph $K$ of this form in
$\bP_+T^*Y$ is a pre-Lagrangian submanifold. The pair
$(K,K)$ has the stable intersection property. This is
a particular case of the Arnold conjecture proved
in [H1],[LS].

\medskip
\noindent
{\bf Example 2.4.B.} Consider the situation described in
1.3.D. Let $L \subset W$ be a Bohr-Sommerfeld Lagrangian
submanifold. Denote by $K$ its full lift to the
prequantization space $QW$. Then $K$ is a pre-Lagrangian
submanifold of $QW$. It is foliated by flat Legendrian lifts
of $L$. Pick up such a lift, say $\Lambda$. If $\pi_2 (W,L) = 0$
then $(K,\Lambda)$ has the stable intersection property.

\medskip
\noindent This example is borrowed from [EHS]
(see also [Ono]) where it is proved that
$(K,\Lambda)$ has the intersection property.
The theory
developed in [EHS] is presented for closed contact manifolds.
In order to apply it to the open manifold
$\Stab M$ and get the
stable intersection property of $(K,\Lambda)$
an additional argument is
needed. We present this argument below.

A prequantization of a symplectic manifold $(W,\Omega)$
is given by the following data:
\begin{itemize}
\item{} a principal $S^1$-bundle
$p: QW \to W$;
\item{} an $S^1$-invariant connection $\xi$ on $QW$
defined by an $S^1$-invariant 1-form
$\alpha$ with
$d\alpha =  p^*\Omega, $ which integrates to $1$ over fibers of
the bundle.
%%%
\end{itemize}
We will denote such a prequantization by
$(QW,p,\alpha)$.
Fix now a sufficiently large positive integer $l$.
Consider the torus $\bT^2 = T^*S^1 / \bZ$
where the action of the group $\bZ$ is generated
by the shift $(r,t) \to (r +l,t)$, and a
symplectic manifold
 $$ W' = (W \times \bT^2, \Omega \oplus dr \wedge dt).$$
Write $B$ for the annulus
$$\{(r,t)\; |\; |r| < l'\},$$
where $0 <l' < l/2$.
We claim that there exists a prequantization
$(QW',p',\alpha')$ of $W'$ with the following properties:
\begin{itemize}
\item{} The bundle $p': QW' \to W'$ over $B$ coincides
with
$$QW \times B \to W \times B,\; (z,r,t) \to  (p(z),r,t).$$
\item{} The form $\alpha'$ over $B$ equals $\alpha + rdt$.
\end{itemize}
Indeed,
denote by $z \to z + s,\; s \in S^1=\bR/\bZ$ the circle action on $QW$.
%%%
The space $QW'$ can be described explicitly as
$$(QW \times T^*S^1)/\bZ,$$
where the $\bZ$-action is generated by the map
\begin{equation}
(z,r,t) \to (z - lt, r+l,t),
\tag{2.4.C}\label{ac}
\end{equation}
and the projection $p'$ is given by
$$p'(z,r,t) = (p(z), r\; {\rm mod}\; l\;,\; t\; {\rm mod}\; 1).$$
One readily checks that $\bZ$-action (2.4.C) preserves
the form $\alpha + rdt$ on $QW \times T^*S^1$.
Define the connection form $\alpha'$ as the push-forward
of $\alpha +rdt$ to
to $QW'$. The claim follows now from the fact that
$QW \times B$ is naturally embedded into the
fundamental domain $QW \times \{-l/2 < r \leq l/2 \}$
of the $\bZ$-action (2.4.C).
\medskip

\noindent
Notice now that open contact manifolds
%%%
$(QW \times B, \alpha + rdt)$ exhaust the stabilization
$\Stab QW$ when $l' \to +\infty$. Thus it suffices
to establish the intersection property of
$(\Stab K, \Stab \Lambda)$ in each such manifold.
In view of then above  claim   $QW \times B$
is contained in $ QW' = Q(W \times \bT^2)$.
Hence, the intersection property for  the pair $(\Stab K, \Stab \Lambda)$ in
 the closed contact manifold $QW'$, which is ensured by
 Theorem 2.5.4 of [EHS], yields   the intersection property
for this pair in $QW \times B$.
%%%
This completes the proof of the stable intersection property.

\section { Relative growth and shape}

\medskip
\noindent
\subsection*{3.1 From paths of contactomorphisms to domains}

\medskip
\noindent
We begin with two lemmas which are proved
at the end of this section (see 3.5 below).

\medskip
\noindent
%\begin{lemma}\label{prop4c}
\proclaim Lemma 3.1.A.
If a  contact diffeomorphism $f\in {\cal D}(M,\xi)$ can be
generated by a positive Hamiltonian, then it can   also be
generated by  a $1$-periodic positive Hamiltonian.
%\end{lemma}

\noindent
Take any element $f\in {\cal D}(M,\xi)$ and represent it by a
path $\{f_t\},\ t\in[0;1],\ f_0=1$ generated by
a 1-periodic contact Hamiltonian $F(x,t)$.  Consider a symplectic
manifold
 $N=SM\times T^*S^1$
endowed with the symplectic form $\omega+dr\wedge dt$ where $(r,t)$,
$r\in\bR,\,t\in\bR/\bZ,$ are canonical coordinates in $T^*S^1$.
Consider domains $$ V^+(\,\{f_t\}\,)=\{r+F(x,t)\geq 0\}\subset N
$$ and $$ V^-(\,\{f_t\}\,)=\{r+F(x,t)\leq 0\}\subset N\enspace. $$

\medskip
\noindent
%\begin{lemma}\label{prop4A}
\proclaim Lemma 3.1.B.
Given two such paths $\{f_t\},\{f_t^\prime\}$ with
$f_1=f^\prime_1$ which are homotopic with fixed endpoints, there
exists a Hamiltonian isotopy of $N$ which takes $V^+(\,\{f_t\}\,)$
to $V^+(\,\{f^\prime_t\}\,)$ and $V^-(\,\{f_t\}\,)$ to
$V^-(\,\{f^\prime_t\}\,)$.
%\end{lemma}

\medskip
\noindent
\subsection*{3.2. Shape}

\medskip
\noindent
Set $$N=S(\mathbb{P}_+T^*\mathbb{T}^n)\times
T^*S^1\subset T^*\mathbb{T}^n\times T^*S^1.$$
 For a  subset $V\subset N$ we define (see \cite{S},\cite{E}) a
  subset $${\rm Shape}(V)\subset H^1(\mathbb{T}^n\times S^1;\,\mathbb{R})
  =H^1(\mathbb{T}^n;\,\mathbb{R})\times \mathbb{R}$$
  as the collection of all pairs $(a,b)$,
   such that there exists a Lagrangian embedding
   $\chi:\mathbb{T}^n\times S^1\to V$ with the
    following properties:
\begin{itemize}
\item[(i)] $\chi^*[pdq+rdt]= (a,b)$;
\item[(ii)] $\chi$ is {\it homologically standard} which means that
the composition of $\chi$ with the natural
projection $V \to \bT^n \times S^1$ induces
 the identity map of
$H^1(\bT^n \times S^1;\bR)$.
\end{itemize}
Let us make identifications
   $$\mathbb{T}^n=\mathbb{R}^n/\mathbb{Z}^n,\; T^*\bT^n = (\bR^n)^* \times \bT^n,
\ H^1(\mathbb{T}^n,\mathbb{R})=
   (\mathbb{R}^n)^*.$$
It was proved by Gromov [Gr1] that every Lagrangian embedding $\chi$ as
above must intersect the split Lagrangian torus $\{p = a,\; r = b\}$.
This result plays a crucial role below.

\medskip
%\begin{example}\label{ex5A}
\noindent{\bf Example 3.2.A
  (Sikorav \cite{S})} Take
a domain ${  U}\subset((\mathbb{R}^n)^*\setminus
\{0\}\,)\times \mathbb{R}$, and set
 $V = U\times \mathbb{T}^n\times S^1$.
  Here $V$ is considered
   as a subset of $N=((\mathbb{R}^n)^*\setminus
\{0\}\,)\times \mathbb{T}^n\times
    \mathbb{R}\times S^1$.
We claim that ${\rm Shape}(V) = U$.
 Indeed, since $V$
is foliated by split Lagrangian tori we get that
$ U \subset {\rm Shape}(V)$. The opposite inclusion is
an immediate consequence of
the abovementioned Lagrangian intersection result.
%    \end{example}
    \medskip

    \noindent Clearly, ${\rm Shape}(V)$ is invariant under Hamiltonian
isotopies of $V$.  Thus, given $f\in {\cal D}$ one can define subsets
$$
\begin{array}{rl}
&{\rm Shape}^+(f)=\,{\rm Shape}(V^+(\,\{f_t\}\,),\\

&{\rm Shape}^-(f)=\,{\rm Shape} (V^-(\,\{f_t\}\,)\,)
\end{array}
$$
which in view of 3.1.B
 do not depend on the particular choice of the path $\{f_t\}$
 representing $f$.
Furthermore, if $V\supset V'$
 then ${\rm Shape}(V)\supset\ {\rm Shape}(V')$.
Combining this with  Proposition 1.4.B
%\ref{prop4B}
  above we get the following
%\begin{proposition}\label{prop5B}

\noindent
\proclaim Proposition 3.2.B. If $f\geq g$ then
${\rm Shape}^+(f)\supset\ {\rm Shape}^+(g)$,
and \\ ${\rm Shape}^-(f)\subset\ {\rm Shape}^-(g)$.
%\end{proposition}

\medskip
\noindent
Note also that ${\rm Shape}^\pm(f)$ is
 invariant under the $\mathbb{R}_+$-action on
  $H^1(\mathbb{T}^n,\mathbb{R})\times\mathbb{R}$.

\medskip\noindent
It is convenient to extract  some numerical invariants
from ${\rm Shape}^\pm(f)$ as follows.  Given $f\in {\cal D}$,
 and $a\in H^1(\mathbb{T}^n,\mathbb{R})\setminus\{0\}$ set
$$
r_-(a,f)=-\inf\{{b}|(a,{b})\in\,{\rm Shape}^+(f)\,\}
$$
and
$$
r_+(a,f)=-\sup\{{b}|(a,{b})\in\,{\rm Shape}^-(f)\,\}\enspace.
$$
Let us list some useful properties of functions $r_+$ and $r_-$.
 Fix $f\in {\cal D}$ and $a\in H^1(\mathbb{T}^n,\mathbb{R})\setminus\{0\}$.

\begin{eqnarray}
(a,{b})\in\,{\rm Shape}^+(f)\ \hbox{for all}\ {b}>-r_-(a,f);\nonumber\\
(a,{b})\in\,{\rm Shape}^-(f)\ \hbox{for all}\ {b}<-r_+(a,f).\label{eqn5C}
\end{eqnarray}

\begin{equation}
r_+(a,f)\geq r_-(a,f)\label{eqn5D}
\end{equation}

\begin{equation}
{\rm (Normalization)}\quad \ r_\pm(a,1)=0\label{eqn5E}.
\end{equation}

\begin{equation}
{\rm (Monotonicity)\quad If}\ f\geq g\ {\rm then}\ r_\pm(a,f)\geq r_\pm(a,g).\label{eqn5F}
\end{equation}

\begin{eqnarray}
&\hbox{(Behaviour under iterations)  For}\ k\in \mathbb{N},\label{eqn5G}\\
&r_-(a,f^k)\geq k r_-(a,f);\nonumber\\
&r_+(a,f^k)\leq kr_+(a,f).\nonumber
\end{eqnarray}

\begin{equation}
r_+(a,f^{-1})=-r_-(a,f)\label{eqn5H}
\end{equation}

\begin{equation}
r_\pm(ca,f)=cr_\pm(a,f)\ \hbox{for all}\ c>0\label{eqn5I}
\end{equation}

\begin{equation}
r_\pm(a,hfh^{-1})=r_\pm(a,f)\   \hbox{for all}\  h\in
{\cal D}\label{eqn5K}
\end{equation}

\medskip
\noindent{\bf Proof of Properties (1) -(8):}

\medskip
\noindent
  Denote by $R_c:N\to N$ the shift
 $r\mapsto r+c$ along the $r$-axis, where $c\in \mathbb{R}$.  Suppose
  that $\chi:\mathbb{T}^n\times S^1\to\{r+F(x,t)\geq 0\}$ is a
homologically standard Lagrangian
  embedding with
  $\chi^*[pdq+rdt]
=(a,b)$.
Then for $c>0\ ,\ R_c\circ\chi(\mathbb{T}^n\times S^1)\subset \{r+F(x,t)>0\}$ and
$(R_c\circ\chi)^*[pdq+rdt]=(a,b+c)$, and hence $(a,{b}+c)\in\,{\rm Shape}^+(f)$.
Taking ${b}$ arbitrary close to $-r_-(a,f)$ we get the first statement
 of (\ref{eqn5C}).  The second one is analogous.
\medskip

\noindent
Recall that a homologically standard Lagrangian torus in $T^*\mathbb{T}^{n+1}$
must intersect a split torus with the same
 Liouville class.
  Since $R_c\circ\chi(\mathbb{T}^n\times S^1)\cap
  \{r+F(x,t)\leq 0\}=\emptyset$ we get that
   $(a,{b}+c)\notin \ {\rm Shape}^-(f)$.  This proves (\ref{eqn5D}).
    The element $1\in{\cal D}$ is represented by $f\equiv 0$.  But
%     (\ref{ex5A})
3.2.A
implies that
      ${\rm Shape}\{r\geq 0\}=H^1(\mathbb{T}^n,\mathbb{R})
      \setminus\{0\}\times\{{b}|{b}\geq 0\}$
       and ${\rm Shape}\{r\leq 0\}=H^1(\mathbb{T}^n,\mathbb{R})
       \setminus\{0\}\times\{{b}|{b}\leq 0\}$.  Thus we get (\ref{eqn5E}).
        The property \ref{eqn5F} is an immediate
         corollary of
%\ref{prop4B}
Proposition 1.4.B.
\medskip

\noindent
Let us show (\ref{eqn5G}).
Assume that $f$ is generated by a 1-periodic Hamiltonian $F(x,t)$.
 Then $f^k$ is generated by $F_k(x,t)=kF(x,kt)$.  Consider the covering
\begin{eqnarray}
\pi:T^*\mathbb{T}^n\times\mathbb{R}
\times\mathbb{R}/\mathbb{Z}&\to& T^*\mathbb{T}^n
\times\mathbb{R}\times\mathbb{R}/\mathbb{Z}\nonumber\\
(x,r,t)&\mapsto& (x,\frac{r}{k}, kt)\enspace.\nonumber
\end{eqnarray}
Clearly, the domain $\{r+F(x,t)\geq 0\}$ lifts to
 $\{r+F_k(x,t)\geq 0\}$.
 If $\chi:\mathbb{T}^n\times S^1\to\{r+F(x,t)\geq 0\}$ is a
homologically standard Lagrangian
 embedding with $\chi^*[pdq +rdt]=(a,b)$
 then $\chi$ lifts to a homologically standard Lagrangian
embedding $\tilde \chi$ with
${\tilde \chi}^*[pdq+rdt] = (a,kb)$.
  Thus if $(a,b)\in\ {\rm Shape}^+(f)$ then
$(a,kb)\in\ {\rm Shape}^+(f^k)$,
  and so $r_-(a,f^k)\geq kr_-(a,f)$.  The second inequality follows
  in the same way.
\medskip

\noindent
The proof of (\ref{eqn5H})
   follows from the fact that $f^{-1}$ is generated by
   $F_{-1}(x,t)=-F(x,-t)$.  The involution
$$
(x,r,t)\mapsto(x,-r,-t)
$$
sends $\{r+F\geq 0\}$ to $\{r+F_{-1}\leq 0\}$, hence the result.
 The property \ref{eqn5I} follows from the fact that contact
 Hamiltonians are homogeneous.  Finally, (\ref{eqn5K}) is obvious.
{\flushright{\hfill \rule{2mm}{2mm}}}
\vfill\eject

%\medskip
%\noindent
\subsection*{3.3 An application to relative growth}

\medskip
\noindent
%\begin{theorem}\label{theorem6A}
\proclaim Theorem 3.3.A. Let $f,g\in{\cal D}$,
   $f$ is a dominant, and $r_-(a,g)>0$ for some
  $a\in H^1(\mathbb{T}^n,\mathbb{R})\setminus\{0\}$.  Then
$$
\gamma(f,g)\geq\frac{r_-(a,g)}{r_+(a,f)}\enspace.
$$
%\end{theorem}

\medskip
\noindent{\bf Proof:}
The proof is based on the properties of functions $r_{\pm}$ listed
in the previous section.
  Set $\gamma_k=\gamma_k(f,g)$ and notice that $\gamma_k>0$.
  Indeed, if $\gamma_k<0$
then since $f^{\gamma_k} \geq g^k$ and $f \geq 1$
we would  have that $g^k \leq 1$, and
$$ 0 \geq r_-(a,g^k) \geq k r_-(a,g) >0, $$
which contradicts to the assumption.
Hence, $f^{\gamma_k}\geq g^k$ implies $$kr_-(a,g)\leq r_-(a,g^k)\leq r_-
 (a,f^{\gamma_k})\leq r_+(a,f^{\gamma_k})\leq \gamma_kr_+(a,f),$$
 and hence $$\frac{\gamma_k}k\geq\frac{r_-(a,g)}{r_+(a,f)}.$$
 Passing to the limit when $k\to+\infty$, we get the desired inequality.
{\flushright{\hfill \rule{2mm}{2mm}}}

\medskip
\noindent{\bf Proof of
%\ref{theorem2B}
Theorem 1.5.A
:}
Identify
$H^1(\mathbb{T}^n,\mathbb{R})$ with
$(\mathbb{R}^n)^*$ where $\mathbb{T}^n=\mathbb{R}^n/
\mathbb{Z}^n$.  It follows from
%\ref{ex5A}
3.2.A
that $r_\pm(p,f)=F(p)$ and $r_\pm(p,g)=G(p)$.  Thus
$\gamma(f,g)\geq\max\limits_{p\neq 0}\frac{G(p)}{F(p)}$
in view of 3.3.A
%\ref{theorem6A}
above.
 The converse inequality follows immediately from the fact that
$$
G(p)\leq\max\frac GF\cdot F(p)
$$
combined with Proposition
%\ref{prop4B}
1.4.B.
{\flushright{\hfill \rule{2mm}{2mm}}}

\medskip
\noindent
\subsection*{3.4 Relative growth, stable norm and minimal action}

\medskip
\noindent
In this section we prove Theorem
%\ref{ex3D}
1.6.E above.
We start from another definition of the stable norm
(see [Gr2, 4.35]).
Let $\rho$ be a Riemannian metric on
  $\mathbb{T}^n$.  For a cohomology
class $a\in H^1(\mathbb{T}^n,\mathbb{R})$ we
   set
   $$\parallel a\parallel^*=\inf\left\{\max\limits_{x\in\mathbb{T}^n}
    |\alpha_x|_{\rho}\ {\Big|}\ \alpha\ \hbox{is a closed 1-form with}\
     [\alpha]=a\right\}.$$
      One can show that $\parallel\ \parallel^*$ is a norm.
       The dual norm $\parallel\ \parallel$ on
$H_1(\mathbb{T}^n,\mathbb{R})$
       is precisely the stable Gromov-Federer norm (see [Gr2]).
        More explicitly,
        for $e\in H_1(\mathbb{T}^n,\mathbb{R}),$
        we have
        $$\parallel e\parallel=\max\{\langle a,e\rangle
        |a\in H^1(\mathbb{T}^n,\mathbb{R}),\ \parallel
a\parallel^*=1\}.$$

\medskip
\noindent
{\bf Proof of Theorem 1.6.E:}
         The geodesic flow on $P_+T^*\mathbb{T}^n$ is given
by the contact
          Hamiltonian  $F(p,q)=$ \ $\parallel p\parallel_\rho$.
Take an element
           $e\in H_1(\mathbb{T}^n,\mathbb{Z})$, and fix $\varepsilon>0$.
            Pick a closed 1-form $\alpha$
            with $[\alpha]=a$ such that
             $\max\limits_{x\in\mathbb{T}^n}
             \parallel \alpha_x\parallel_\rho\leq 1+\varepsilon$
              and $\langle a,e\rangle=\parallel e\parallel$.
               Since ${\rm graph}(\alpha)\subset \{F\leq 1+\varepsilon\}$
               we have
               ${\rm graph}(\alpha)\times \{r=-1-\varepsilon\}\subset
               \{r+F\leq 0\}$.  Thus
$(a,-1-\varepsilon)\in\ {\rm Shape}^-(f)$,
from which we conclude that $r_+(a,f)\leq 1+\varepsilon$.

\medskip
\noindent Recall now that $H_1(\mathbb{T}^n,\mathbb{Z})$ is identified
with a subgroup of $\pi_1(\Gamma(M,\xi)\,)$.  With this identification
and also identifying $(\mathbb{R}^n)^*$ with $H^1(\mathbb{T}^n,\mathbb{R})$
 we get that $e$ is generated by a Hamiltonian $G(p)=\langle
 p,e\rangle$.
  In view of  Example 3.2.A
% \ref{ex5A}
we have
   $r_\pm(a,e)=\langle a,e\rangle =\parallel e\parallel$, and  applying
   Theorem 3.3.A
% \ref{theorem6A}
we get that
$$
\gamma(f,e)\geq\frac{r_-(a,e)}{r_+(a,f)}\geq\frac{\parallel e\parallel}{1+\varepsilon}\enspace.
$$
Since $\varepsilon$ is arbitrarily small, this implies
 the desired inequality.
{\flushright{\hfill \rule{2mm}{2mm}}}

\medskip
\noindent
It was established by Bangert see
( \cite{Ba}, Sect. 2.C) that the stable norm is related
to Mather's minimal action [Ma].  Namely, assume that $\mu$ is an
invariant Borel probability measure of the geodesic flow on $T\bT^n$.
 Define its action $$A(\mu)=\frac 12 \int\limits_{TM} |v|^2_\rho d\mu(v).$$
   Denote by $\cM$ the set of all measures with $\cA(\mu)<\infty$.
    The rotation number $R(\mu)$ is an element from $H_1(\bT^n,\bR)$
     which satisfies
      $\langle\,[\alpha],R(\mu)\,\rangle=\int\limits_{TM}\alpha d\mu$,
       where a closed 1-form is considered as a function on $TM$.
       Define Mather's {\it minimal action}
$$
\beta:H_1(\bT^n,\bR)\to\bR
$$
by the formula
$$\beta(e)=\inf\{A(\mu)\,{\Big |}\,{\mu\in \cM}\;
{\rm and}\; {R(\mu)= e} \}.$$
  Bangert proved that $\beta(e)=\frac 12  \parallel e\parallel^2_\rho$.
   Thus our result above implies that
$\gamma(f,e)\geq \sqrt{2\beta(e)}$.
This inequality
is not a specific feature of geodesic flows. In fact
it remains true for any $f \in \cD(\bP_+T^*\bT^n)$ generated
by a positive {\it time-independent} contact Hamiltonian $F(p,q)$
whose square $F^2$ is strictly convex with respect to $p$-variable.
This is an easy consequence of [CIPP, Cor. 1].
It is unclear however whether an inequality of this type
remains true if one considers time-one maps of {\it time-dependent}
Hamiltonians. Further there is a strong feeling that
the quantity $\gamma (f,e)$ is related to
invariant pre-Lagrangian tori of $f$
(compare with Siburg's theory [Si] which links together invariant
tori, Mather's minimal action and Hofer's geometry).

%\vfill\eject
\medskip
\noindent
\subsection*{3.5 Proof of lemmas 3.1.A and 3.1.B}

\medskip
\noindent
{\bf Proof of Lemma 3.1.A:}
Let $f_t,\,t\in[0,1],$ be a path of
contactomorphisms connecting $1$ with $f$, and ${\tilde f}_t:SM\to SM$ the
symplectization of this path. Let $F(x,t),\,t\in [0,1],\,x\in SM$,
be a positive homogeneous Hamiltonian which generates the path
${\tilde f}_t$. We will assume that $F(x,t)$ is extended for  $t\in [0,2]$
as a positive Hamiltonian,  set $G(x,t)=F(x,t+1),\,t\in[0,1],$ and
denote by $g_t:M\to M$ and ${\tilde g}_t:SM\to SM$ the contact and
symplectic isotopies defined by the Hamiltonian $G(x,t)$. Take the
product $S_2M=SM\times SM$ with the symplectic structure
$\Omega=(-\omega)\oplus\omega$ and consider Lagrangian graphs
$$\Gamma_t=\{(x,{\tilde f}_t(x))|\,x\in SM\}\quad\hbox{and}\quad
\Delta_t=\{(x,{\tilde g}_t(x))|\,x\in SM\}$$  of   symplectomorphisms $f_t$
and $g_t$. Notice that these graphs are invariant with
respect to the diagonal
action of $\bR_+$ on
 $(S_2M,\Omega)$. Furthermore, the
$\bR_+$-homogeneous Hamiltonian functions $\widetilde F(x,y,t)=F(y,t)$
 and $\widetilde G(x,y,t)=G(y,t)$, $(x,y)\in S_2M,
 \,t\in[0,1]$,
 generate symplectic isotopies
${\rm Id}\times {\tilde f}_t:SM\to SM$ and ${\rm Id}\times {\tilde g}_t:SM\to SM$
  which move  the
 diagonal $\Gamma_0=\Delta_0$ to
$\Gamma_t$ and $\Delta_t$, respectively.
In particular, we see that   the
 Lagrangian isotopies $\Gamma_t$ and $\Delta_t$ are
generated by   positive
 (homogeneous) Hamiltonian functions. The converse is also true: if
 $B_t$ is a homogeneous Lagrangian
graphical (with respect to the splitting
 $S_2M=SM\times SM$) isotopy, which  is
 generated by a family of positive homogeneous  functions
 $K_t:\Gamma_t\to\bR$ then the corresponding family of
 $\bR_+$-equivariant symplectomorphisms of $SM$ is generated by
 a positive homogeneous (time-dependent) Hamiltonian on $SM$.

\medskip
\noindent Next,  observe  that there exists a $\bR_+$-invariant
neighborhood $U$ of the diagonal
  $\Gamma_0$ which is equivariantly
symplectomorphic to a $\bR_+$-invariant neighborhood of the $0$-section
  in the cotangent bundle $T^*(SM)$. Here we canonically  extend
  the action of $\bR_+$ on $SM$ to a
conformally symplectic action on $T^*(SM)$ as
follows. Write $R_c$ for the action $q \to cq$,
where $q \in SM$ and $c \in \bR_+$. Then the extended action
in the canonical coordinates $(p,q)$ on $T^*(SM)$ is given by
$$(p,q) \to (c(R_{c}^*)^{-1}p, R_{c}(q)).$$
   There exists  a small $\delta>0$ such that for  $t\leq \delta$
   the Lagrangian submanifolds
    $\Gamma_t$ and $\Delta_t$ are contained in $U$, and hence
     we can identify
    them with    Lagrangian
submanifolds of $T^*(SM)$,  still denoted by $\Gamma_t$
     and $\Delta_t$. These submanifolds are
    $\bR_+$-invariant, and hence   are defined by
    homogeneous generating functions $S_t:SM\to\bR$ and $T_t:SM\to\bR$,
respectively.
     Notice that the Lagrangian
isotopies $\Gamma_t$ and $\Delta_t$, $t\in[0,\delta]$
   are    generated by positive  Hamiltonian functions
    $$A_t(p,q)=\frac{\partial S_t}{\partial t}(q)\quad\hbox{and}\quad
    B_t(p,q)=\frac{\partial T_t}{\partial t}(q), $$ and hence we
    have $\frac{\partial S_t}
{\partial t}(q),\frac{\partial T_t}{\partial t}(q) >0$
    for all $q\in
    SM,\,t\in[0,\delta]$.
Let us recall that $\Gamma_0$ and $\Delta_0$
(viewed as submanifolds of the cotangent bundle of $SM$)
 coincide with the
zero-section in $T^*(SM)$.
Thus $S_0 \equiv T_0 \equiv 0$. The inequalities above imply
that for $t \in (0,\delta]$ the functions $S_t$ and $T_t$
are positive. Further,
for a sufficiently   small positive
$\varepsilon<\frac{\delta}2$ we have $T_\varepsilon
<S_{\frac{\delta}2}$. Therefore, one can find  a smooth  family of
$\bR_+$-homogeneous functions $U_t,\,t\in[0,\delta],$ on $SM$ which
coincides with $T_t$ for $t\in[0,\varepsilon]$ and with $S_t$ for
$t\in[\frac{\delta}2 ,\delta]$, and such that $\frac{\partial
U_t}{\partial t}>0$ for all $t\in[0,1]$. If $\delta $ have been chosen
  small enough  then one can guarantee, in addition, that the
Lagrangian submanifolds  of $T^*(SM)$ generated by the functions
$U_t$, and viewed in a neighborhood of the
diagonal in $S_2M$, are graphical
with respect to the splitting $S_2M=SM\times SM$,
and hence correspond to  $\bR_+$-equivariant
symplectomorphisms ${\tilde h}_t:SM\to SM$.   It remains to observe that
the family ${\tilde h}_t,\,t\in[0,\delta],$  together with the family
${\tilde f}_t,\,t\in[\delta,1]$, define a smooth path of
$\bR_+$-equivariant symplectomorphisms $SM\to SM$ which is
generated by a time-periodic positive
Hamiltonian. {\hfill \rule{2mm}{2mm}}
\bigskip

\medskip
\noindent
{\bf Proof of Lemma 3.1.B:}
 Let $h_{t,s}$ be a family of homogeneous symplectomorphisms of $SM$ such that
\begin{itemize}
\item[$*$] for  a fixed $s$,  the symplectomorphism
$h_{t,s}$ is generated by a 1-periodic Hamiltonian $F(x,t,s)$;
\item[$*$] $h_{0,s}\equiv 1\ ,\ h_{1,s}\equiv f$;
\item[$*$] for a fixed $t,\ h_{t,s}$ is generated by $G(x,t,s)$.
\end{itemize}
  We assume without loss of generality that
$h_{t+1,s}=h_{t,s}\circ f$ for all $s$, and hence  $G$ is 1-periodic in $t$.
  We also have  $\frac{\partial F}{\partial s}=\frac{\partial G}{\partial t}-\{F,G\}$.  Consider $G(x,t,s)$ as
  a Hamiltonian on $SM\times T^*S^1(x,r,t)$, depending on time $s$.  Thus
   $$\frac{dr}{ds}=-\frac{\partial G}{\partial t},\ \frac{dt}{ds}=0,\ \frac{dx}{ds}={\rm sgrad} G_{t,s}.$$
    We claim that this Hamiltonian flow takes $\{r+F(x,t,0)=0\}$ to $\{r+F(x,t,s)=0\}$,
    in other words $r(s)+F(x(s),t,s)=0$.  Indeed, differentiating by $s$ one gets
$$
-\frac {\partial G}{\partial t}+\{F,G\}+\frac{\partial F}{\partial
s}=0\enspace.
$$
In fact, the Hamiltonian isotopy in question is given explicitly by
$$
(x,r,t)\mapsto \left(h_{t,s}h^{-1}_{t,0}x,\ r-\int^s_0\frac{\partial G_{t,u}}{\partial t}\left(h_{t,u}(h^{-1}_{t,0}(x))\right)du,\ t\right)\enspace.
$$
Note that it is homogeneous with respect to $\mathbb{R}_+$-action on
$
SM\times T^*S^1$:  $$ (x,r,t)\mapsto(cx,cr,t).
$$
{\flushright{\hfill \rule{2mm}{2mm}}}

\medskip
\noindent
{\bf Acknowledgments.} We thank Miguel Abreu and Ana Cannas da Silva,
the organizers
of the Symplectic Geometry conference at
Lisboa in June 1999 for the opportunity to present results of
this paper. We are grateful to Alexander Givental for illuminating
discussions concerning the non-linear Maslov index.
A part of the manuscript has been written while
the second author was visiting
IHES (Bures-sur-Yvette). He thanks IHES for the hospitality.

\end{document}